\documentclass[a4paper,12pt,reqno]{amsart}

\pdfoutput=1

\usepackage[margin=2.5cm, headsep=1cm, footskip=1cm]{geometry}
\usepackage[hang, flushmargin]{footmisc}
\usepackage[english]{isodate}
\usepackage{amsmath, amsthm, amssymb, mathtools}
\usepackage{enumitem}
\usepackage{xcolor}
\usepackage{tikz, tikzsymbols, tikz-cd}
\usepackage{multicol}
\usepackage[colorlinks=true, urlcolor=black, citecolor=black, linkcolor=black, hyperfootnotes=true]{hyperref}
\usepackage[capitalise, nameinlink, noabbrev, nosort]{cleveref}
\usepackage{wasysym} % for alchemy >:3 
\usepackage[normalem]{ulem} % for crossing out 
\usepackage{stmaryrd} 

\usepackage{matlab-prettifier} % for code 
\usepackage{appendix} % for appendices 

\usepackage{float}
\usepackage{comment}

\newcommand{\Mod}[1]{\ (\mathrm{mod}\ #1)}
\definecolor{indigo}{HTML}{492DA5}

\allowdisplaybreaks[2]
\setenumerate{listparindent=\parindent}

\providecommand{\noopsort}[1]{}

\makeatletter\g@addto@macro\bfseries{\boldmath}\makeatother

\makeatletter
\let\origsection\section
\renewcommand\section{\@ifstar{\starsection}{\nostarsection}}
\newcommand\sectionspace{\vspace{0.5ex}}
\newcommand\nostarsection[1]{\sectionspace\origsection{#1}\sectionspace}
\newcommand\starsection[1]{\sectionspace\origsection*{#1}\sectionspace}
\makeatother

\crefname{page}{page}{pages}

\setlist[enumerate]{font=\normalfont}
\crefname{enumi}{}{}
\crefname{enumii}{}{}

\numberwithin{equation}{section}
\crefname{equation}{Equation}{equations}

\crefname{condition}{condition}{conditions}
\creflabelformat{condition}{#2(#1)#3}

\newtheorem{theorem}{Theorem}[section]

\crefname{thm}{Theorem}{Theorems}

\newtheorem{lemma}[theorem]{Lemma}
\crefname{lemma}{Lemma}{Lemmas}

\newtheorem{prop}[theorem]{Proposition}
\crefname{prop}{Proposition}{Propositions}

\newtheorem{cor}[theorem]{Corollary}
\crefname{cor}{Corollary}{Corollaries}

\theoremstyle{definition}

\newtheorem{definition}[theorem]{Definition}
\crefname{definition}{Definition}{Definitions}

\newtheorem{remark}[theorem]{Remark}
\crefname{remark}{Remark}{Remarks}

\crefname{remarks}{Remarks}{Remarks}

\newtheorem{example}[theorem]{Example}
\crefname{example}{Example}{Examples}

\newcommand{\Z}{\mathbb{Z}}
\newcommand{\Q}{\mathbb{Q}}
\newcommand{\C}{\mathbb{C}}
\newcommand{\F}{\mathbb{F}}
\newcommand{\Zmodd}{\mathbb{Z}/2\mathbb{Z}}

\newcommand{\Zmodp}{\mathbb{Z}/p\mathbb{Z}}
\newcommand{\Zmodn}{\mathbb{Z}/n\mathbb{Z}}

\newcommand{\GGG}{\mathcal{G}}
\newcommand{\SSS}{\mathcal{S}}

\DeclareMathOperator{\supp}{supp} 

\newcommand{\boundellipse}[3]% center, xdim, ydim
{(#1) ellipse (#2 and #3)
}

\begin{document} 

\date{\today}
\title[Singular ideals of the cyclic-headed snakes]{Singular ideals over arbitrary fields for the cyclic-headed snakes} 

\author[Abdellatif]{Ramla Abdellatif}
\author[Clark]{Lisa Orloff Clark}
\author[Jansen]{Roy Jansen}
\author[Marsland]{Stephen Marsland}


\thanks{This research collaboration was initiated by the first two named authors at the 2023 International Algebra Conference in the Philippines (IACP). Many thanks to Prof. Jocelyn P. Vilela and the other organizers of this conference for providing us this opportunity to meet in such great conditions. Thanks also to VUW and UPJV for facilitating our visits to each other, and to the anonymous reviewer for very useful feedback.}

\address[L.O. Clark, R. Jansen, S. Marsland]{School of Mathematics and Statistics, Victoria University of Wellington, PO Box 600, Wellington 6140, NEW ZEALAND}
\email{\href{mailto:lisa.orloffclark@vuw.ac.nz}{lisa.orloffclark@vuw.ac.nz}}
\email{\href{mailto:jansenroy@myvuw.ac.nz}{roy.jansen@vuw.ac.nz}}
\email{\href{mailto:Stephen.Marsland@vuw.ac.nz}{stephen.marsland@vuw.ac.nz}}
\address[R. Abdellatif]{Laboratoire Ami\'{e}nois de Math\'{e}matique Fondamentale et Appliqu\'{e}e, Universit\'{e} de Picardie Jules Verne, 33 rue Saint-Leu - 80 039 Amiens Cedex 1, FRANCE}
\email{\href{mailto:ramla.abdellatif@math.cnrs.fr}{ramla.abdellatif@math.cnrs.fr}}

\begin{abstract} 
    We study the Steinberg algebras with coefficients in an arbitrary field $K$
    for the cyclic-headed snake groupoids, which are basic examples of non-Hausdorff groupoids. We are particularly interested in elements of this algebra that are no longer continuous, known as singular functions. These functions form an ideal, which may contain proper subsets that are themselves ideals of the Steinberg algebra. We provide three conditions %on the choice of the base field $K$ or the number of heads of the snake,% 
    under which such `proper subset' ideals exist: first, when the number of heads of the snake divides the characteristic of the base field; second, when the base field is of non-prime characteristic; and third, when certain cyclotomic polynomials split over the base field. We also show the existence of many further subset ideals not covered by these conditions.  

    We fully explore the cases of the two- and three-headed snakes. In the three-headed snake, we prove that the ideal of singular functions properly contains non-zero ideals of the Steinberg algebra if, and only if,  the base field $K$ is a splitting field of $\Phi_3(T) = T^{2}+T+1$, the third cyclotomic polynomial. Consequently, there are always proper subset ideals when $K$ has characteristic a prime not congruent to $-1 \Mod 3$. 
\end{abstract}

\maketitle 

%\tableofcontents 

\section{Introduction}

Steinberg algebras are constructed from certain topological groupoids and have proven useful in purely algebraic settings as models for inverse semigroup algebras \cite{Steinberg2010}, Leavitt path algebras \cite{LPAbook}, or Kumjian–Pask algebras \cite{CP}. Steinberg algebras were introduced independently in \cite{Steinberg2010} and \cite{CFST}. Complex Steinberg algebras -- that is, with coefficients in $\mathbb{C}$ -- have also found applications in functional analysis (see for example \cite{BCFS} or, more recently, \cite{BGHL}); in particular, they sit densely inside the associated groupoid C*-algebra \cite{snake}, and the ability to move between algebraic and analytic frameworks has become a valuable technique for making progress in both areas.

At first, the specific choice of $\C$ as field of coefficients did not appear to be particularly significant in the algebraic setting. However, when the underlying groupoid is not Hausdorff, the study of the simplicity of such algebras \cite{Simplicity, Nek, StSz, StSz2} made it clear that the choice of the field of coefficients, and in particular of its characteristic, plays a crucial role. 

Given a field $K$ and a topological groupoid $\mathcal{G}$, recall that the Steinberg algebra $A_{K}(\mathcal{G})$ over $K$ is defined as the $K$-linear span of characteristic functions of compact open bisections of $\mathcal{G}$.  In a non-Hausdorff groupoid, we have compact sets that are not closed; and as a result, the associated Steinberg algebra may contain functions whose support has empty interior: such functions are called \textit{singular functions}. From \cite[Theorem~3.14]{Simplicity}, we know that a Steinberg algebra that contains singular functions cannot be simple, as the set of singular functions is a proper ideal. We are here interested in whether $A_{K}(\GGG)$ is \emph{S-simple} -- that is, whether its ideal of singular functions has non-trivial proper subsets that are also ideals of $A_{K}(\GGG)$. We call such ideals \emph{proper subset ideals}. The S-simplicity of $A_{K}(\GGG)$ is a very different concept to the simplicity of the ideal of singular functions itself. In the snake groupoids considered here, the ideal of singular functions is in fact never simple for snakes with more than two heads. 

In this paper, we begin to study the relationship between the field $K$ and the S-simplicity of $A_K(\mathcal{G})$ when $\mathcal{G}$ is non-Hausdorff. We consider the Steinberg algebras of the groupoids known as the $n$-headed snakes, $\mathcal{G}_{n S}$, which we review in \cref{prelims}.  

The two-headed snake $\GGG_{2S}$ \cite[Example~2.1]{snake} is one of the most basic examples of a non-Hausdorff groupoid. For us, the body of the snake comprises the well-known Hausdorff Cantor set, which is the unit space of this groupoid. The fact that the groupoid $\GGG_{2S}$ is non-Hausdorff is concentrated at one location: the two heads of the snake, denoted as $0$ and $1$,  which interact with each other as in $\mathbb{Z} / 2 \mathbb{Z} $. The open sets of $\GGG_{2S}$ are the open sets of the Cantor set together with copies of those sets with head $0$ replaced by head $1$. When we add more heads to the snake, which interact with each other as in a finite group $G$, the fact that the corresponding $G$-headed snake $\GGG_{G S}$ is non-Hausdorff is again confined to the heads. 

In \cref{algs} we show that the singular functions of the Steinberg algebra $A_{K}(\GGG_{G S})$ are precisely those whose support is contained in $G \subseteq \mathcal{G}_{G S}$ -- that is, functions that are non-zero only on the heads of the snake. We also note the distinction between ideals of the Steinberg algebra and subideals of the ideal of singular functions, and compare with the group ring $KG$. In \cref{cyclicheaded_section} we restrict to the cyclic-headed snakes, where the group $G$ of heads is a cyclic group $\Zmodn$; the groupoid $\mathcal{G}_{n S}$ is then the $n$-headed snake. 
In the two-headed snake, we show that $A_K(\mathcal{G}_{2 S})$ is always S-simple, regardless of the field $K$. However, with three or more heads, the S-simplicity of $A_K(\mathcal{G}_{n S})$ is field-dependent. In \cref{ThreeConditionsSection}, we provide three conditions under which $A_K(\mathcal{G}_{n S})$ is not S-simple: when $\mathrm{char}(K) | n$; when $n$ is not prime; and when $K$ is a \emph{splitting field} for the $n$th cyclotomic polynomial $\Phi_n(T)$. 
We also show that these three conditions do not provide an exhaustive list, and we lay out a computational method for identifying further subset ideals. We provide an example in the five-headed snake of a subset ideal not covered by our conditions. 

For the three-headed snake algebra $A_{K}(\mathcal{G}_{3 S})$, we \emph{are} able to provide a complete list of all possible subset ideals, in \cref{Examples: Two and Three Headed Snakes}. We prove that $\mathcal{S}_{K}(\mathcal{G}_{3S})$ is not S-simple if and only if $K$ is a splitting field for $\Phi_3(T)$ -- that is, $K$ contains the roots of $T^2+T+1$. We conclude by providing a classification of such fields. 

In future work, we hope to investigate criteria for S-simplicity in more complicated non-Hausdorff groupoids, not necessarily of snake type.

\section{Preliminaries} \label{prelims}
\subsection{Ample groupoids and their Steinberg algebras} 
We provide some brief preliminaries for context. For a more detailed presentation of the relevant groupoids and their Steinberg algebras, the interested reader can refer to, for example, \cite{Steinberg2010}.  

A \textit{groupoid} $\mathcal{G}$ is a generalisation of a group where the binary operation is only partially defined. As a result, the identity, or `unit', is not a single element, but a set of elements.  We let $\mathcal{G}^{(0)}$ denote the set of all units in $\mathcal{G}$.  Each $\gamma \in \mathcal{G}$ has an inverse $\gamma^{-1} \in \mathcal{G}$ and we write $\textbf{r},\textbf{s}:\mathcal{G} \to \mathcal{G}^{(0)}$ for the maps
\[\textbf{r}(\gamma)=\gamma\gamma^{-1} \quad \text{ and } \quad \textbf{s}(\gamma)=\gamma^{-1}\gamma \ , \]
which we call the \emph{range} and \emph{source} maps respectively.  Composition between two elements $\gamma$ and $\eta$ in $\mathcal{G}$ is defined if, and only if,  $\textbf{s}(\gamma) = \textbf{r}(\eta)$. 

A motivating example for this paper comes from bundles of groups: given a non-empty set $I$ and, for each $i \in I$, a group $\Gamma_{i}$, 
\[ 
\mathcal{G} := \bigsqcup_{i \in I}\Gamma_i 
\] 
is a groupoid. A pair of elements of $\mathcal{G}$ is composable if and only if there is some $i \in I$ such that both elements are in $\Gamma_{i}$. Composition is then defined as composition in the group $\Gamma_{i}$. 

A \emph{topological groupoid} is a groupoid equipped with a topology such that composition and inversion define continuous maps. An \emph{open bisection} is an open set $B \in \mathcal{G}$ such that the source and range maps on $B$ restrict to homeomorphisms. An \emph{ample} groupoid has a topology generated by a basis of compact open bisections. 

Steinberg algebras are constructed from certain topological groupoids that generalise the class of discrete groups. The unit space of our groupoids must be Hausdorff, but the main point of our investigation is that the groupoid itself need not be. 

Given a field $K$ and an ample groupoid $\mathcal{G}$ with Hausdorff unit space, the \textit{Steinberg algebra} of $\mathcal{G}$ over $K$ is the $K$-algebra $A_{K}(\mathcal{G})$ of functions $f : \mathcal{G} \rightarrow K$ of the form:  
    \[ f = \sum\limits_{B \in F} a_{B} 1_{B} \ , \] 
    where $F$ is a finite collection of compact open bisections $B$ of $\mathcal{G}$, each $a_{B} \in K$, and $1_{B}$ denotes the characteristic function of $B$. 
    Addition and scalar multiplication in $A_{K}(\mathcal{G})$ are defined pointwise, with convolution given by: 
    \begin{equation} \label{convformgeneral}
        \forall \  f, g \in A_{K}(\GGG), \ \ \gamma \in \GGG, \ \ \ \ \ (f*g)(\gamma) := \sum\limits_{\mathbf{r}(\eta)= \mathbf{r} (\gamma)} f(\eta) g(\eta ^{-1} \gamma) \ .
    \end{equation}
When there is no ambiguity which groupoid $\mathcal{G}$ we are working with, we write $A_{K}$ for the algebra $A_{K}(\GGG)$. 

\subsection{The Two-Headed Snake} 
The two-headed snake is a well-known example of a non-Hausdorff space. We study a totally disconnected version that has a groupoid structure, as in \cite[Example~2.1]{snake}.
Let $X$ be the Cantor set, viewed as a subset of $[0,1] \in \mathbb{R}$:  this serves as the body of our snake. 

\begin{remark}
    One could as easily define the snake using a different body $X$, such as $X = \{0\} \cup \{ 1/n: n>0 \}$ with the appropriate (ample, Hausdorff) topology. Since our investigation focuses on the heads of the snake (being the location of any unusual non-Hausdorff behaviour), we can disregard most details of the body. For convenience, we use the Cantor set definition of the snake given in \cite{snake}. 
\end{remark}

The topology on $X$ has a basis $\mathcal{B}_{X}$ consisting of cylinder sets $B$, each of which is compact in $X$. $X$ is Hausdorff and compact with this topology. We let the point $0_{\mathbb{R}} = 0$ be the `head' of the snake, and add another head denoted by $1$. Note that this additional head $1$ is not the point $1_{\mathbb{R}}$ at the other end of the snake. 

\begin{figure}[H]
\centering
\begin{tikzpicture}[node distance={15mm}, thick, main/.style = {draw, circle}, every edge quotes/.append style = {font=\footnotesize}]

\draw[line width=2pt] (0,0) to (7,0) ; 
\draw (0,0) node[circle,fill,label=left:{0},label=below:{$0_{\mathbb{R}}$}] {}; 
\draw (0,2) node[circle,fill,label=left:{1}] {}; 
\draw (7,0) node[circle,fill,label=below:{$1_{\mathbb{R}}$}] {}; 
\draw[dashed, line width=2pt] (0,2) to [bend right, looseness=1] (1,0); 

\draw \boundellipse{0,1}{1}{2};

\draw (3.5,-0.25) node[label=below:{$X$}] {};

\end{tikzpicture} 
\caption{ Schematic of the two-headed snake groupoid} 
\end{figure}
Thus, as a set, the two-headed snake groupoid is $\mathcal{G}_{2 S} = X \cup \{ 1 \}.$ The composable pairs in $\GGG_{2S}$ are given by:
  \[
        \mathcal{G}_{2 S}^{(2)} = \{ (\gamma,\gamma) : \gamma \in X \} \cup  \{ (0,1) , (1,0), (1,1) \} \ . 
  \]
Composition is denoted $+$ and defined by: 
\[ \left\{
\begin{array}{l} 
\gamma +\gamma = \gamma \ \text{ for all } \gamma \in X,\\  
0+1 =1 \ ,\\ 
1+0 =1 \ , \\
1+1 =0 \ , 
\end{array} \right. \] 
so that the heads $\{0,1\}$ interact as in $\Zmodd$. 
The basis of the topology on $\mathcal{G}_{2 S}$ consists of the basic open sets $\mathcal{B}_X$ of the Cantor set; plus copies of each $B \in \mathcal{B}_{X}$ that contains the point $0$, with $0$ replaced by the additional head $1$. That is:
    \begin{equation*}
        \mathcal{B}_{2 S} = \mathcal{B}_X \cup \{ (B \backslash \{ 0 \} ) \cup \{ 1 \} : B \in \mathcal{B}_X \text{ such that } 0 \in B \}.
    \end{equation*}
This groupoid can be viewed as a group bundle consisting of trivial groups at every point in the body $X$ and one copy of $(\mathbb{Z} / 2 \mathbb{Z}, +)$ at the head.  
We use subscripts to denote whether a set contains the head 0 or the head 1, where relevant. In particular, $X_{0} = X$ is the snake with head 0 (and not 1) and $X_{1}$ is the snake with head 1 (and not 0). 

With this structure, $\mathcal{G}_{2 S}$ is an ample groupoid with Hausdorff unit space. The groupoid is not Hausdorff because the elements 0 and 1 cannot be separated by open sets.  Each element of the basis $\mathcal{B}_{2 S}$ is a compact open bisection, while any set that contains more than one head of the snake cannot be a bisection; since, as $\textbf{r}(1) = \textbf{r}(0) = 0$, the range map is not injective on such sets (and likewise for the source map). Note that sets of the form $X_{a} \backslash B_{a}$, where $a \in \Zmodd$, are compact open bisections that are not closed.

\subsection{The Group-Headed Snakes} \label{Section_NonHausdorffSnakes}

A natural way to generalize the two-headed snake is to add more heads. Let $G$ be a finite group: this will be the group of heads. For convenience we write $G$ additively, with $id_G = 0_G$, although we do not require $G$ to be Abelian. We define the $G$-headed snake groupoid $\mathcal{G}_{G S}$ in a similar manner as the two-headed snake. 

\begin{figure}[ht]
\centering
\begin{tikzpicture}[node distance={15mm}, thick, main/.style = {draw, circle}, every edge quotes/.append style = {font=\footnotesize}]

\draw[line width=2pt] (0,0) to (7,0) ; 
\draw (0,0) node[circle,fill,label=left:{$0_G$},label=below:{$0_{\mathbb{R}}$}] {}; 
\draw (0,1) node[circle,fill,label=left:{$j$}] {}; 
\draw (0,2) node[circle, fill,label=left:{$j'$}] {}; 
\draw (0,2.8) node[label=left:{}] {}; 
\draw (7,0) node[circle,fill,label=below:{$1_{\mathbb{R}}$}] {}; 

\draw[dashed, line width=2pt] (0,1) to [bend right, looseness=0.7] (1.5,0); 
\draw[dashed, line width=2pt] (0,2) to [bend right, looseness=0.5] (1.5,0); 
\draw[dashed, line width=2pt] (0.1,2.8) to [bend right, looseness=0.4] (1.5,0); 

\draw[line width=3pt, line cap=round, dash pattern=on 0pt off 2\pgflinewidth, line width=5pt] (0,2.5) to (0,3.3); 

\draw \boundellipse{0,1.3}{1.7}{2.5}; 

\draw (3.5,-0.25) node[label=below:{$X$}] {}; 

\end{tikzpicture} 
\caption{Schematic of the group-headed snake groupoid} 
\end{figure}

As a set, the group-headed snake is $\mathcal{G}_{G S} = X \cup G$, where $0 \in X$ matches with $0_{G} \in G$, and the composable pairs are: 
\[
    \mathcal{G}_{G S}^{(2)} = \{ (\gamma,\gamma) : \gamma \in X \} \cup ( G \times G ) . 
\]
Composition is given by $\gamma + \gamma = \gamma$ for all $\gamma \in X$ and by the group composition for all points in $G$. Generalising the two-headed snake, the basis for the topology is: 
\[
\mathcal{B}_{G S} = \mathcal{B}_X \cup \bigsqcup\limits_{j \in G} \{ (B \backslash \{ 0 \} ) \cup \{ j \} : B \in \mathcal{B}_X \text{ such that } 0 \in B \} . 
\] 
As in the two-headed snake, we denote by $X_j$ the snake with head $j$ and no other heads. Again, $\GGG_{GS}$ is an ample, non-Hausdorff groupoid, and any set that contains more than one head of the snake is not a bisection. Likewise, for any head $j \in G$, the set $X_{j} \backslash Z(\mu)_{j}$ is a compact open bisection that is not closed.

\section{Steinberg Algebras of the Snakes} \label{algs} 

Let $G$ be a finite group and $K$ a field. Consider the Steinberg algebra $A_K(\mathcal{G}_{G S})$ of the $G$-headed snake. Let $f \in A_{K}$. We denote by $f_j \in K$ the value of function $f$ on head $j \in G \subseteq \mathcal{G}_{G S}$; that is, $f_j := f(j)$. 

\begin{lemma} \label{formoffunc} 
Let $f \in A_{K}(\mathcal{G}_{G S})$. Then there exists a finite collection $F_{w}$ of compact open bisections of $\GGG_{GS}$ that contain no heads of the snake, and a collection $\{f_{j} : j \in G \} \cup \{f_{B} : B \in F_{w}\}$ of scalars in $K$, such that $f$ can be written in the form: 
\begin{equation*}  
    f = \displaystyle \sum_{j \in G} f_{j} 1_{X_{j}} + \sum\limits_{B \in F_{w}} f_{B} 1_{B} \ .
\end{equation*} 
\end{lemma} 

\begin{proof}
Begin by writing $f$ in the standard form: $f = \sum\limits_{B \in F} a_{B} 1_{B}$, where $F$ is a finite collection of compact open bisections of $\GGG_{GS}$. Without loss of generality, we can assume each $B \in F$ is a basic set (by decomposing $B$ into a disjoint union of basic open sets and writing $a_{B}1_{B}$ as the corresponding sum). 

Let $|G| = n$. Split $F$ into $n+1$ disjoint parts: the collection of those $B$ that contain the head $0$ (denoted $F_0$); those that contain the head $1$ (denoted $F_1$); and so forth for all heads $j \in G$; and those $B$ that contain no heads at all (denoted $F_w$). 
Then: 
\[
    f = \sum\limits_{j \in G} \left( \sum\limits_{B \in F_{j}} a_{B}1_{B} \right)  + \sum\limits_{B \in F_{w}} a_{B} 1_{B} \ . 
\] 

Consider any $F_{j}$. For each $B \in F_{j}$, note that $j \in B$ (so that $B = B_j$) and thus $X_{j} \backslash B$ is a compact open bisection. Write: 
\[ a_{B} 1_{B} = a_{B} 1_{X_{j}} - a_{B}  1_{X_{j} \backslash B} \ ,\] 
and so: 
\[
    \sum\limits_{B \in F_{j}} a_{B} 1_{B} = \left(\sum\limits_{B \in F_{j}} a_{B}\right) 1_{X_{j}} + \left(\sum\limits_{B \in F_{j}} - a_{B}  1_{X_{j} \backslash B}\right).
\]
Note that the head $j \notin X_j \backslash B$, and so move $\sum\limits_{B \in F_{j}} - a_{B}  1_{X_{j} \backslash B}$ into the sum over $F_{w}$, with $F_{w}$ adjusted accordingly. Now $F_j = \{X_j\}$ only, and so the coefficient $\sum\limits_{B \in F_{j}} a_{B}$ must be equal to $f(j) = f_{j}$. 

Repeating this process for all $j \in G$, adjusting $F_{w}$ accordingly, we have $f = \sum\limits_{j \in G} f_{j} 1_{X_{j}} + \sum\limits_{B \in F_{w}} a_{B} 1_{B}$ as claimed. 
\end{proof} 

As $G$ is finite, we can fix an ordering of its $n$ elements. Given $f_{0},f_{1}, ..., f_{n-1} \in K$, we let: 
\[ [f_{0},f_{1},...,f_{n-1}] := f_{0} 1_{X_{j_0}} + f_{1} 1_{X_{j_1}} + ...+ f_{n-1} 1_{X_{j_{n-1}}} \ . \]
Then by \cref{formoffunc} any $f \in A_K(\GGG_{G S})$ can be written as: 
\[ f = [f_{0},f_{1},...,f_{n-1}] + \sum\limits_{B \in F_{w}} f_{B} 1_{B} \ . \] 

\begin{example}
In the two-headed snake, for $f,g \in A_K(\GGG_{2 S})$, the convolution formula \cref{convformgeneral} simplifies to: 
    \begin{align*}
        (f*g)(\gamma) = \begin{cases}
            f(\gamma) \cdot g(\gamma) \text{ if } \gamma \in X \backslash \{0\} \\ 
            f_0g_0 + f_1g_1\text{ if } \gamma = 0 \\ 
            f_0g_1 + f_1g_0 \text{ if } \gamma = 1 \ . 
        \end{cases} \\ 
    \end{align*} 
\end{example} 
More generally, in the $G$-headed snake, for $f,g \in A_K(\GGG_{G S})$: 
    \begin{align} \label{convform_Gh}
        (f*g)(\gamma) = \begin{cases}
            f(\gamma) \cdot g(\gamma) \text{ if } \gamma \in X \backslash \{0\} \\ 
            \sum\limits_{j \in G} f_j g_{j^{-1} i} \text{ if } \gamma = i \in G \ . 
        \end{cases} 
    \end{align}     
In the sequel, we are interested in the values on the heads of the snake. Where relevant, we simply write $\sum_{B \in F_{w}} a_{B} 1_{B}$ for the values of $(f*g)$ on the body.  
Note that when the group of heads $G$ is Abelian, convolution in $A_K(\mathcal{G}_{GS})$ is commutative.

\subsection{Ideals of Singular Functions of the Snakes} \label{Singular Ideals Intro Section}

When a groupoid $\mathcal{G}$ is not Hausdorff, the Steinberg algebra $A_K(\mathcal{G})$ may contain some functions that are not continuous, called \emph{singular functions}. 
For a function $f:\mathcal{G} \to K$, the \emph{support} of $f$ is: 
\[
\supp(f)\coloneqq \{\gamma \in \mathcal{G} : f(\gamma) \neq 0\}.
\]
This is sometimes referred to as the `open support' when the function $f$ is continuous. However, if $f$ is not continuous, then $\supp(f)$ is not necessarily open.  A function $f \in A_K(\mathcal{G})$ is \emph{singular} if $\supp(f)$ has empty interior. The set of all singular functions in $A_K(\mathcal{G})$ is:
    \begin{equation*}
    \mathcal{S}_{K}(\mathcal{G}) := \{ f \in A_K(\mathcal{G}) : \supp(f) \text{ has empty interior} \}.
    \end{equation*}  
We know from \cite[Proposition~3.7]{Simplicity} that $\mathcal{S}_{K}(\mathcal{G})$ is an ideal in $A_K(\mathcal{G})$. 
For the snakes, singular functions are precisely the functions in $A_K$ that have no support on the body of the snake. 
\begin{lemma} \label{2h_singnosuppbody}
    A function $f \in A_K(\mathcal{G}_{GS})$ is singular if and only if $\supp(f) \subseteq G$. 
\end{lemma} 

\begin{proof}
    Let $f \in A_K(\mathcal{G}_{G S})$. Assume there exists $\gamma \in X \backslash \{0\}$ in $\supp(f)$. Write $f$ in the form of \cref{formoffunc}. Then there exists some compact open bisection $U \in F_w$ such that $\gamma \in U$. 
    Let $F_{\gamma} \subseteq F_w$ be the subcollection of $B \in F_w$ such that $\gamma \in B$. Then: 
    \[ 
        V := U \cap \left( \bigcap\limits_{B \in F_{\gamma}} B \; \; \bigg\backslash  \bigcap\limits_{B \in F_w \backslash F_{\gamma}} B \right) 
    \] 
    is an open set, and $V \subseteq \supp(f)$. Thus $\supp(f)$ does not have empty interior, so $f \notin \SSS_{K}(\GGG_{G S})$.  
    Finally, note that $G \subseteq \mathcal{G}_{GS}$ contains no open sets, so if $\supp(f) \subseteq G$ then $f$ is by definition singular. 
\end{proof} 

\cref{2h_singnosuppbody} tells us that our singular functions can be written $f = [f_{0},f_{1},...,f_{n-1}]$, with the condition that $f$ vanishes on the body of the snake:  
\begin{lemma} \label{singular_headssumtozero} 
    The function $f = [f_{0},f_{1},...,f_{n-1}] \in A_K(\mathcal{G}_{G S})$ is singular if and only if $f_{0}+f_{1}+...+f_{n-1} = 0_{K}$.
\end{lemma} 

\begin{proof}
    Choose $\gamma \in X \backslash \{0\}$. Then $\gamma \in X_j$ for all $j \in G$, so $f(\gamma) = f_{0}+f_{1}+...+f_{n-1}$. By \cref{2h_singnosuppbody}, the set of singular functions is precisely the functions such that $f(\gamma) = 0_K$ for all $\gamma \in X \backslash \{0\}$.  
\end{proof}
\cref{singular_headssumtozero} says: 
\[
\mathcal{S}_K(\mathcal{G}_{G S}) = \{ [f_{0},f_{1},...,f_{n-1}] : f_{0}+f_{1}+...+f_{n-1} = 0_{K} \} . 
\] 

Note that $\mathcal{S}_K$ is always non-empty and non-zero -- for example, for any $k \in K$ we have $[k,-k,0,0,...,0] \in \mathcal{S}_K$ -- and so $A_K (\mathcal{G}_{G S})$ is always non-simple.

\subsection{Comparison with the group ring $KG$}

In this section, we show that the ideal of singular functions is isomorphic to an ideal of the group ring $KG$ called the \textit{augmentation ideal}. Group rings are well known (see for example \cite{DummitandFoote}) and provide the structure needed to prove \cref{allprincipal}. 

Let $(G, \cdot)$ be a group and $K$ a field. Recall that the group ring $KG$ consists of the elements $\{ f = \sum_{j \in G} f_j j: f_j \in K \}$, with addition and scalar multiplication pointwise. The product between $f,g \in KG$ is defined as: 
    \[
        fg = \sum\limits_{i \in G} \left( \sum\limits_{j \in G} (f_j g_{j^{-1} \cdot i}) i \right) . 
    \] 
When $G$ is Abelian, $KG$ is called the group algebra over $K$.
Consider:  
\[
    C := \left\{ \sum_{j \in G} f_j 1_{X_j}: f_j \in K \right\} \subseteq A_K( \mathcal{G}_{G S} ) , 
\] 
 the functions that are constant on the body of the snake. Note that $C$ is a subalgebra of $A_K$. We claim that $C$ is isomorphic, algebraically (that is, without topology), to the group ring  $K G$, similar to \cite[Section 3]{StSz}. To see this, recall that for compact open $U,V$ we have $1_U * 1_V = 1_{UV}$ where $UV = \{u v : u \in U \mathbf{s} (V), v \in \mathbf{r} (U) V \}$. In particular, for any $i,j \in G$, we have $1_{X_i} * 1_{X_{j}} = 1_{X_{i \cdot j}}$. The map $j \mapsto 1_{X_j}$ is thus a homomorphism between $(G, \cdot)$ and $\left( C, * \right)$, and thence an isomorphism. From \cref{convform_Gh}, observe that for functions $f = \sum_{j \in G} f_j 1_{X_j}$ and $g = \sum_{j \in G} g_j 1_{X_j}$ in $C$, we have: 
    \[
        (f*g) = \sum\limits_{i \in G} \left( \sum\limits_{j \in G} (f_j g_{j^{-1} \cdot i}) 1_{X_i} \right) \ . 
    \]
Thus we have an isomorphism $\sum_{j \in G} f_j j  \longmapsto \sum_{j \in G} f_j 1_{X_j}$, as required. 
\begin{definition}
The augmentation map $\varepsilon: K G \rightarrow K$ is $\sum_{i \in G} f_{i} i \mapsto \sum_{i \in G} f_{i}$, and the kernel of $\varepsilon$ is known as the \emph{augmentation ideal} of $KG$.  
\end{definition}
\begin{lemma} \label{SK_is_augmentideal} 
    The ideal of singular functions of the group-headed snake $\mathcal{S}_K(\mathcal{G}_{G S})$ is algebraically isomorphic to the augmentation ideal of $K G$. 
\end{lemma} 

\begin{proof}
    By \cref{singular_headssumtozero}, the ideal $\mathcal{S}_K$ is precisely the functions $f = \sum_{j \in G} f_j 1_{X_j}$ such that $\sum_{j \in G} f_j = 0_K$; that is, (isomorphic to) the kernel of $\varepsilon$. 
\end{proof} 
The ideals of $\mathcal{S}_K$ are thus subideals of the augmentation ideal.

\subsection{Subideals of $A_K$ vs. Ideals of $\mathcal{S}_K$} \label{subidvssubsetid} 

However, an ideal of $\mathcal{S}_K$ is not necessarily also an ideal of $A_K$. Studying the group ring $KG$ alone is not sufficient to determine the structure of the Steinberg algebra $A_K(\mathcal{G}_{G S})$.  
\begin{definition}
    We call an ideal $J \unlhd A_K$ such that $J \subseteq \mathcal{S}_K$ a \emph{subset ideal} of $A_K$; and if $J \subsetneq \mathcal{S}_K$, we say that $J$ is a \emph{proper} subset ideal of $A_K$. 
    If $A_K$ has no proper subset ideals, we say that $A_K$ is \emph{S-simple}. 
\end{definition} 
We must be careful to distinguish between ideals of $A_K$ and ideals of $\mathcal{S}_K$. We use the notation $\langle f \rangle_{A_K}$ to denote the ideal of $A_K$ generated by $f$, and likewise $\langle f \rangle_{\mathcal{S}_K}$ the ideal of $\mathcal{S}_K$. Note that as $\mathcal{S}_K \subseteq A_K$, we have for all $f \in A_K$: 
\[
    \langle f \rangle_{\mathcal{S}_K} \subseteq \langle f \rangle_{A_K} \ . 
\]
Our focus in this paper is on S-simplicity and subset ideals of $A_K$, as opposed to the simplicity and subideals of $\mathcal{S}_K$. In \cref{nonsimplicity_of_SK}, for the cyclic $n$-headed snakes with $n > 2$, we show an ideal of $\mathcal{S}_K$ that is not an ideal of $A_K$, and demonstrate that $\mathcal{S}_K$ is never simple for these snakes. 

The collection of subset ideals of $A_K$ is contained in, but not always equal to, the collection of subideals of $\mathcal{S}_K$. 
In general, if $S$ is an ideal of a ring $R$ and if $I$ is a subideal of $S$, then $I$ is not necessarily an ideal of $R$. That is: 
\[
    I \unlhd \mathcal{S}_K \unlhd A_K \nRightarrow I \unlhd A_K. 
\]  
Thus an ideal of $\mathcal{S}_K$ is not necessarily also a subset ideal of $A_K$. So if $\mathcal{S}_K$ is not simple, then $A_K$ may or may not be S-simple. 
On the other hand, a subset ideal $J$ of $A_K$ is always an ideal of $\mathcal{S}_K$: 
\[
    J \unlhd A_K \text{ and } J \subseteq \mathcal{S}_K \Rightarrow J \unlhd \mathcal{S}_K \ .  
\] 
To see this, observe that as $J \subseteq \mathcal{S}_K$ is an ideal (of $A_K$), $J$ must be closed under addition and scalar multiplication. Choose $f \in \mathcal{S}_K$; then $f \in A_K$ also, so $f*g \in J$ for all $g \in J$, and so $J \unlhd \mathcal{S}_K$. 
\begin{remark} \label{subsetid_isof_KG}
     Similarly, $J$ is also an ideal of the subalgebra $C$: because $C \subseteq A_K$, all $f \in C$ are in $A_K$, and so $f*g \in J$ for all $g \in J$. Thus, by the isomorphism of \cref{SK_is_augmentideal}, all subset ideals $J$ are (isomorphic to) ideals of the group ring $KG$. Again, the reverse does not hold; not every ideal of the group ring is (isomorphic to) a subset ideal $J$.  
\end{remark} 
The group ring alone does not give us enough structure to identify which ideals of $\mathcal{S}_K$ are also subset ideals of $A_K$. However, when $K G$ is well-known, it can provide a useful starting point. For example, in \cref{ThreeConditionsSection}, we explore the cyclic-headed snakes, where $G = \Zmodn$. The group algebra $K \Zmodn$ is isomorphic to $K [T] \backslash \langle T^n - 1 \rangle$, and so has ideals corresponding to any roots of $T^n - 1$ that are in $K$. In \cref{nheads_cycloideals}, we find that such ideals -- when they exist -- are in fact also proper subset ideals of $A_K (\mathcal{G}_{n G})$.

\section{Cyclic-headed Snakes} \label{cyclicheaded_section} 

From this point we deal only with the cyclic-headed snakes. Let the group of heads $G = \Zmodn$ be the cyclic group of order $n$, written additively, and consider the cyclic $n$-headed snake $\mathcal{G}_{nS}$. Using \cref{convform_Gh} we write: 
\begin{equation} \label{convform_nh}
    (f*g) = \left[ \sum\limits_{j \in G} f_j g_{n-j}, \sum\limits_{j \in G} f_j g_{1-j}, ..., \sum\limits_{j \in G} f_j g_{n-1-j} \right] + \sum\limits_{B \in F_{w}} a_{B} 1_{B} \ . 
\end{equation}

\begin{prop} \label{allprincipal}
    All subset ideals of the $n$-headed snakes are principal. 
\end{prop} 

\begin{proof}
    When $G = \Zmodn$ a cyclic group, it is well-known (see for example \cite{DummitandFoote}) that the group algebra $K \Zmodn$ is a principal ideal domain regardless of the characteristic of $K$. By \cref{subsetid_isof_KG}, all subset ideals of $A_K$ are also ideals of $C \cong KG$, and are therefore singly generated by an element of $C$. As $C \subseteq A_K$, subset ideals are also principal in $A_K$. 
\end{proof}

The ideal of singular functions $\mathcal{S}_K$ is also principal. The augmentation ideal is known to be generated by $\{ j - 1_{\Zmodn} : j \in \Zmodn \} \cong \{ 1_{X_j} - 1_{X_1}: j \in \Zmodn \}$ -- see for example \cite[Theorem 12.4]{Johnson}. As we are interested in ideals of $A_K$, and not in all of the subideals of the augmentation ideal, we use a different proof which for our purposes is more concise. 
\begin{theorem}
\label{fullideal}
    The ideal $\mathcal{S}_K(\mathcal{G}_{n S}) \unlhd A_K (\mathcal{G}_{n S})$ is generated by $[1,-1,0,...,0]$: 
    \[ \mathcal{S}_K(\mathcal{G}_{n S}) = \langle [1,-1,0,...,0] \rangle_{A_K} . \]
\end{theorem} 

\begin{proof} 
    As $1-1+0+...+0 = 0$ regardless of field, $[1,-1,0,...,0]$ is singular, and so $\langle [1,-1,0,...,0] \rangle_{A_K} \subseteq \mathcal{S}_K$. 
    For the reverse inclusion, choose any $s \in \mathcal{S}_K$. We show that there exists some $f \in A_K$ such that $s = f*[1,-1,0,...,0]$ -- in fact, there are arbitrarily many such $f$. 
    Now $f = [ s_0, \ \  s_0+s_1, \ \  s_0+s_1+s_2, \ \  ..., \ \ s_0+...+s_{n-2}, \ \ 0] \in A_K$. Since $s \in \mathcal{S}_K$, by \cref{singular_headssumtozero} we have $s_0+...+s_{n-2} = -s_{n-1}$. Using \cref{convform_nh} we calculate: 
    \begin{align*}
        f*[1,-1,0,...,0] =& [f_0 - f_{n-1}, \ \ f_1 - f_0, \ \ f_2 - f_1, \ \ ..., \ \ f_{n-1}-f_{n-2}] \\ 
        =& [s_0 - 0, \ \ (s_0+s_1) - s_0, \ \ ..., \ \ 0 - (s_0+...+s_{n_2})] \\ 
        =& [s_0, s_1, ... , s_{n-1}] = s 
\end{align*}  
    as needed. Thus $s \in \langle [1,-1,0,...,0] \rangle_{A_K}$ as required. 
\end{proof} 

The following is immediate. 
\begin{cor} \label{propersubid_oneonenotin} 
    Let $J$ be a subset ideal of $A_K$, that is, $J \unlhd A_K$ and $J \subseteq \mathcal{S}_K$. Then $J$ is a proper subset ideal, that is, $J \neq \mathcal{S}_K$, if and only if $[1,-1,0,...,0] \notin J$. 
\end{cor} 
Thus to check whether a subset ideal $J \subseteq \mathcal{S}_K$ is a \emph{proper} subset of $\mathcal{S}_K$, it suffices to check whether there exists $s \in J$ and $f \in A_K$ such that $s*f = [1,-1,0,...,0]$. 
The choice of $[1,-1,0,...,0]$ as the generator of $\mathcal{S}_K$ is convenient but certainly not unique -- see \cref{licitvariations} and \cref{otherSKgenerators_remark}.

\subsection{Licit variations of the generating function} \label{licitvariations}

Consider a singular function $s = [s_0, s_1, ..., s_{n-1}] \in \mathcal{S}_K(\mathcal{G}_{n S})$. When we are working with $\langle s \rangle_{A_K}$ we can make adjustments to our choice of $s$ without altering the ideal. 
Let $1_{X_j} = [0,...,0,1,0,...,0] \in A_K(\mathcal{G}_{nS})$, where the only non-zero term is at head $j$. 
For all $f \in A_K$, convolution with $1_{X_j}$ has the effect of cyclically permuting the values on the heads of $f$ along $j$-many steps to the left. We omit the proof as it is a straightforward computation. 
\begin{lemma} \label{permutationallowed}
    Let $f \in A_K(\mathcal{G}_{n S})$ and $j \in \Zmodn$. Then the $(k-j)$th term of $f$ becomes the $k$th term of $1_{X_j} * f$: 
    \[ (1_{X_j} * f)_k = f_{k-j} \] 
    for all $k \in \mathbb{Z}_n$. 
\end{lemma} 
Observe that for all $j \in \Zmodn$, the function $1_{X_j}$ is in $A_K$ (and not in $\mathcal{S}_K$). 
Thus, when considering ideals of $A_K$, it does not change the ideal if the coefficients of any generating function are rotated, so long as their relative positions are preserved. For example, 
\[
    \langle [1,-1,0,...,0] \rangle_{A_K} = \langle [0,1,-1,0...,0] \rangle_{A_K} \ . 
\] 
For nonzero $s \in \mathcal{S}_K$ we can therefore, without changing any ideal generated by $s$, require that $s_0 \neq 0$. 
As all ideals are closed under scalar multiplication, we can further require that $s_0 = 1$, and adjust the other terms of $s$ by a factor of $\frac{1}{s_0}$. 
Finally, by \cref{singular_headssumtozero} we know that the heads of $s$ must sum to zero. We can thus without loss of generality always write: 
\begin{equation} \label{licitvariations_EQUATION}
    s = [1,s_1, ..., s_{n-2}, -(1+s_1+...+s_{n-2})] \ . 
\end{equation}
Note that it is in general not possible to otherwise permute the terms of $s$ without altering $\langle s \rangle_{A_K}$.

\subsection{Non-Simplicity of $\mathcal{S}_K$} \label{nonsimplicity_of_SK}

Our focus in this paper is on whether subset ideals of $A_K$ exist, not an investigation of subideals of $\mathcal{S}_K$ itself. Indeed, for the cyclic-headed snakes with more than two heads, there \emph{are} such subideals, regardless of $K$ and of $\Zmodn$. 

\begin{lemma} \label{SKneversimple}
    Over any field $K$, for any $n > 2$, the subalgebra $\mathcal{S}_K (\mathcal{G}_{n S})$ is not simple. 
\end{lemma}

\begin{proof}
    Consider $\langle [1,-1,0,...,0] \rangle_{\mathcal{S}_K}$. We show that this is always a proper subideal of $\mathcal{S}_K$. 
    As in the proof of \cref{fullideal}, by \cref{convform_nh} we have for all $f \in \mathcal{S}_K$: 
    \[
        f * [1,-1,0,...,0] = [f_0 - f_{n-1}, f_1 -f_0, ..., f_{n-1} - f_{n-2}] \ . 
    \] 
    Therefore, a function $s$ is an element of $\langle [1,-1,0,...,0] \rangle_{\mathcal{S}_K}$ if and only if there exists some $f \in \mathcal{S}_K$ such that $s = [f_0 - f_{n-1}, f_1 -f_0, ..., f_{n-1} - f_{n-2}]$. Recall from \cref{singular_headssumtozero} that $f \in \mathcal{S}_K$ must satisfy $f_0+f_1+...+f_{n-1} = 0_K$. 

    However, the functions $[0,1,-1,0,...,0]$ and $[1,0,...,0,-1]$ are always in $\mathcal{S}_K$, regardless of $K$. We claim that at most one of them is contained in $\langle [1,-1,0,...,0] \rangle_{\mathcal{S}_K}$. 
    
    First, if  $[0,1,-1,0,...,0] \in \langle [1,-1,0,...,0] \rangle_{\mathcal{S}_K}$, then there exists $f \in \mathcal{S}_K$ such that: 
    $f_1 - f_0 = 1$, and $f_2 - f_1 = -1$, and $f_j = f_{j-1}$ for all $j \in [3,n]$. Beginning with $j = n = 0$ and ending with $j = 3$ we have: 
    \[
        f_{0} = f_{n-1} = f_{n - 2} = ... = f_{3} = f_{2}  
    \] 
    and with $f_1 = f_0 +1$, all required equations hold. 
    Hence $f = [f_0, f_0 +1, f_0, ..., f_0]$, for any $f_0 \in K$. For this $f$ to be singular, we must have $nf_0 + 1 = 0_K$; that is, $\mathrm{char}(K) \equiv 1\Mod n$. 

    Similarly, if $[1,0,...,0,-1] \in \langle [1,-1,0,...,0] \rangle_{\mathcal{S}_K}$, then we must have $f \in \mathcal{S}_K$ such that: $f_0 - f_{n-1} = 1$, and $f_{n-1} - f_{n-2} = -1$ and $f_j = f_{j-1}$ for all $j \in [1, n-2]$. 
    So we must have $f = [f_0, f_0, ..., f_0 - 1]$ for any $f_0 \in K$; and for $f \in \mathcal{S}_K$ we must have $n f_0 - 1 = 0_K$, and thus $\mathrm{char}(K) \equiv -1\Mod n$. 

    When $n > 2$, $\mathrm{char}(K)$ cannot be congruent to both $1$ and $-1\Mod n$, and so we cannot have both $[0,1,-1,0,...,0]$ and $[1,0,...,0,-1]$ in $\langle [1,-1,0,...,0] \rangle_{\mathcal{S}_K}$. Therefore $\langle [1,-1,0,...,0] \rangle_{\mathcal{S}_K} \subsetneq \mathcal{S}_K$, and so $\mathcal{S}_K$ is always non-simple.  
\end{proof} 

\begin{remark} \label{2h_alwaysssimple}
    The two-headed snake algebra $A_K(\mathcal{G}_{2 S})$, when $n = 2$, is the `trivial' case, where there are neither proper subset ideals of $A_K$ nor proper subideals of $\mathcal{S}_K$. To see this, let $J \subseteq \mathcal{S}_{K}(\mathcal{G}_{2 S})$ be a nonzero ideal of $A_K(\mathcal{G}_{2 S})$; recall from \cref{subidvssubsetid} that $J$ is thus any ideal of $\mathcal{S}_K$. Choose a nonzero $s \in J$. Then as $s$ is a singular function, by \cref{singular_headssumtozero} we must have $s = [s_0, -s_0] = s_0 [1,-1]$ with $s_0 \neq 0$. Then $\frac{1}{s_0} \in K$, so $[1,-1] \in J$, and so $J = \mathcal{S}_K$. 
\end{remark}

From now we deal only with ideals of $A_K$.

\section{Some Conditions for Non-S-Simplicity of $\mathcal{S}_K (\mathcal{G}_{n S})$} \label{ThreeConditionsSection}

We now present three forms of proper subset ideals, that occur under certain conditions on the field $K$ and the group of heads $\Zmodn$. \cref{summarycharts} provides a quick summary of the ideals and when they occur. In \cref{FurtherSubsetIdeals}, we show that these three kinds of ideal are not an exhaustive list of possible proper subset ideals. 

Our methodology is similar to the proof of \cref{SKneversimple}. Consider any $s \in \mathcal{S}_K$. By \cref{propersubid_oneonenotin} we have that $\langle s \rangle_{A_K} \subsetneq \mathcal{S}_K$ if and only if there is no $f \in A_K$ such that $f * s = [1,-1,0,...,0]$. We use \cref{convform_nh} to calculate $f*s$ for a general $f \in A_K$, and compare terms with $[1,-1,0,...,0]$, to obtain conditions on $f$ in terms of $s$. For certain choices of $s$ this method is straightforward to do by hand. In \cref{FurtherSubsetIdeals} we explore a generalised approach.

\subsection{Field characteristic dividing $n$}

The subset ideal introduced here is in fact always an ideal of $A_K (\mathcal{G}_{n S})$. It is only a subset of $\mathcal{S}_K$ when $\mathrm{char}(K) | n$. 

\begin{prop} \label{oneoneone_subideal}
    If $\mathrm{char}(K) | n$, then $A_K (\mathcal{G}_{n S})$ has a proper subset ideal given by: 
    \[
    \langle [1,1,....,1] \rangle_{A_K} \subsetneq \mathcal{S}_K(\mathcal{G}_{n S}) \ . 
    \]
\end{prop}

\begin{proof}
    Choose any $f \in A_K(\mathcal{G}_{n S})$. Then by \cref{convform_nh}: 
    \begin{align*}
        f*[1,1,....,1] =& [f_0+f_1+...+f_{p-1}, ..., f_0+f_1+...+f_{n-1}] \\ 
        =& \sum\limits_{j \in \Zmodn} f_j [1,1,...,1]  \ .
    \end{align*} 
    And so: 
    \begin{equation} \label{oneall_scalarmultiples}
        \langle [1,1,....,1] \rangle_{A_K} = \{ k \cdot [1,1,...,1] : k \in K \} \ .
    \end{equation}
    Regardless of the choice of field, there is no $k \in K$ such that $k \cdot [1,1,...,1]$ equals $[1,-1,0,...,0]$. Therefore by \cref{propersubid_oneonenotin} we have $\langle [1,1,1,....,1] \rangle_{A_K} \neq \mathcal{S}_K (\mathcal{G}_{n S})$. 
    Finally, by \cref{singular_headssumtozero}, the function $[1,1,....,1]$ is singular if and only if $1+1+1+...+1 = n = 0_K$; that is, $\langle [1,1,....,1] \rangle_{A_K} \subsetneq \mathcal{S}_K(\mathcal{G}_{n S})$ if and only if $\mathrm{char}(K) | n$. 
\end{proof}

\subsection{Non-prime heads}

The proper subgroups of cyclic group $\Zmodn$ are each isomorphic to $\mathbb{Z} / j \mathbb{Z}$ where $j$ is a proper divisor of $n$. When such proper subgroups exist (that is, when $n$ is not prime), there are associated proper subset ideals of $A_K (\mathcal{G}_{n S})$. 

\begin{prop} \label{oneminusone_ideal}
    If $n$ is not prime, then $A_K (\mathcal{G}_{n S})$ has proper subset ideals generated by $1_{X_0} - 1_{X_j}$, for each $j$ which nontrivially divides $n$, that is: 
    \[
        \langle [1,0,...,0,-1,0,...,0] \rangle_{A_K} \subsetneq \mathcal{S}_K (\mathcal{G}_{n S}) \ . 
    \]
\end{prop} 

\begin{proof}
    Let $j$ be a nontrivial divisor of $n$; that is, $\langle j \rangle = \{ 0,j,2j,...,n-j \} $ is a proper subgroup of $\Zmodn$. Choose any $f \in A_K(\mathcal{G}_{nS})$. Then by \cref{convform_nh}:  
    \begin{align*}
        f*(1_{X_0} - 1_{X_j}) =& [f_0 - f_{n-j}, f_1 - f_{n-j+1}, ..., f_{n-1} - f_{n-j-1}] \ . 
    \end{align*} 
    Assume for a contradiction that $\langle 1_{X_0} - 1_{X_j} \rangle_{A_K} = \mathcal{S}_K$. Then by \cref{propersubid_oneonenotin} there exists some $f$ such that $f * (1_{X_0} - 1_{X_j}) = [1,-1,0,...,0]$. We must have $f_0 - f_{n-j} = 1$, and $f_{k} = f_{k-j}$ for all $k \in [2, n-1]$. Beginning with $k = n-j$ and ending with $k = j$ we have: 
    \[
        f_{n-j} = f_{n-2j} = f_{n - 3j} = ... = f_{2j} = f_{j}  
    \] 
    and as $ j \in [2,n-1]$, we also have $f_{j} = f_{0}$, so $f_{n_j} = f_{0}$; a contradiction. 
\end{proof}

\begin{remark} \label{otherSKgenerators_remark}
    If $j$ is any generator of $\Zmodn$ (in particular, if $j \nmid n$), then 
    \[ 
        \langle 1_{X_0} - 1_{X_j} \rangle_{A_K} = \langle [1,0,...,0,-1,0,...,0] \rangle_{A_K} = \mathcal{S}_K \ , 
    \]
    by a similar proof as of \cref{fullideal}. When $n = p$ a prime, every $j \neq 0 \in \Zmodp$ is a generator. 
\end{remark}

\subsection{Splitting field of $\Phi_q(T)$ for $q|n$} \label{splittingfield_subsec}

We now introduce some elementary number theory that we use in the sequel, and provide a brief overview of some standard results. For more detail, see for example \cite{DummitandFoote}. 

\subsubsection{Cyclotomic Polynomials}

\begin{definition} \label{cyclopol_def}
    Given a field $K$, the $n$th roots of unity in $K$ are the roots of the equation $x^n-1=0$ in $K$. An $n$th root of unity is `primitive' if it is not also an $m$th root of unity for any $m < n$. 
\end{definition}
Typically, $\zeta_n$ denotes a primitive $n$th root of unity, and $\omega_n$ any (primitive or not) $n$th root of unity. Every $n$th root of unity $\omega_n$ is a primitive $d$th root of unity $\zeta_d$ for some $d | n$.  
Note that if $\zeta_n$ is a primitive $n$th root of unity, then every other $n$th root of unity is given by $\omega_n = \zeta_n^k$ for some $k \le n$ -- in particular, if a field $K$ contains one primitive root $\zeta_n$, then $K$ contains all $n$th roots of unity. If $k$ is coprime to $n$, then $\zeta_n^k$ is also a primitive $n$th root of unity. 

\begin{definition}
    Let $K$ be a field containing the primitive $n$th roots of unity $\zeta_n$. 
    The $n$th cyclotomic polynomial over $K$ is the unique irreducible polynomial with integer coefficients whose roots are precisely every $\zeta_n \in K$: 
    \[ \Phi_n(T) = \prod_{\substack{1 \le k < n \\ k \nmid n}} (T - \zeta_n^k) \ . \] 
    Equivalently:  
    \[ \Phi_n(T) = \frac{T^n-1}{\prod_{\substack{1 \le q < n \\ q \mid n}} \Phi_q(T) } \ . \]
\end{definition}
For prime $p$, the $p$th cyclotomic polynomial is $\Phi_p(T) = T^{p-1}+...+T^2+T+1$. The number of primitive $n$th roots of unity is given by Euler's totient function $\varphi(n)$, which for a prime $p$ is $\varphi(p) = p-1$. Observe that when $\mathrm{char}(K) = p$, we have that $1$ is a root of $\Phi_p(T) \in K[T]$. 

\begin{definition}
    A field $K$ is called a \emph{rupture field} for a polynomial $f(T)$ if at least one of the roots of $f(T)$ is an element of $K$. If $K$ contains all the roots of $f(T)$, then $K$ is a \emph{splitting field} for $f(T)$. 
\end{definition} 
If a field $K$ is a rupture field for a cyclotomic polynomial $\Phi_n(T)$, then $K$ is always a splitting field for $\Phi_n(T)$. 
A field $K$ is said to be \emph{algebraically closed} when every non-constant polynomial with coefficients in $K$ has a root in $K$. Thus the algebraic closure of any field is, by definition, a splitting field for all $\Phi_n(T)$. 

We define the polynomial 
    \[
        \Psi_n(T) := T^{n-1} + ... + T^2 + T + 1 \ . 
    \] 
As $T^n - 1 = (T-1)(T^{n-1}+...+T+1)$ and $\Phi_1(T) = T-1$, observe that:  
    \[ 
        \Psi_n(T) = \prod_{\substack{q \neq 1 \\ q \mid n}} \Phi_q(T) \ . 
    \]
That is, the roots of $\Psi_n(T) \in K[T]$ are precisely the primitive $q$th roots of unity $\zeta_q$ for all $q \mid n$, $q \neq 1$; including the root $\zeta_q = 1$ when $\mathrm{char}(K) | n$.

\subsubsection{Cyclotomic-form subset ideals} 

In $A_K(\mathcal{G}_{n S})$, there are proper subset ideals associated with any roots of $\Psi_n(T)$ -- that is, any $\zeta_q$ with $q>1$, $q |n$ -- that are in $K$. Let $q,d \in \mathbb{N}$ such that $q > 1$ and $q | d | n$. Then we have $m \in \mathbb{N}$ such that $n = md$. For each primitive $q$th root of unity $\zeta_q$, we define: 
\[
    c(\zeta_q, d) := [ 1,0,...,0, \ \ \overbrace{ \zeta_q,0,...,0,}^{m \text{ terms }} \ \ \zeta_q^2,0,...,0, \ \ ..., \ \ \zeta_q^{d-1},0,...,0 ] \ . 
\] 
That is,  
\[
    c(\zeta_q, d)_i = \begin{cases}
                        \zeta_q^k \text{ if } i = km \text{ for } k \in [0,d-1] \\ 
                        0 \text{ otherwise. }
                    \end{cases} 
\] 

\begin{comment}
\begin{remark}
    Recall from \cref{SK_is_augmentideal} that $\mathcal{S}_K$ is isomorphic to the augmentation ideal of the group algebra $K\Zmodn$. Therefore all subideals of $mathcal{S}_K$ are also subideals of $K \Zmodn$. 
    
    It is well known that $K \Zmodn \cong K[T] \backslash \langle T^n - 1 \rangle$. With $T^n - 1 = (T-1) \Psi_n(T)$, the subideals of $\langle T^n - 1 \rangle$ in $K [T]$ -- and thus the subideals of $K \Zmodn$ -- are determined by the divisors of $\Psi_n(T)$.

    , we expect there to be subideals of $\mathcal{S}_K$ corresponding with divisors of $\Psi_n(T)$. In this case, these subideals of $\mathcal{S}_K$ are also subset ideals of $A_K$. 
\end{remark} 
\end{comment} 

\begin{theorem} \label{nheads_cycloideals} 
    Consider the $n$-headed snake. If $K$ is a splitting field for $\Phi_q(T)$ for any $q > 1$ such that $q | n$, then $A_K (\mathcal{G}_{nS})$ has proper subset ideals, given by:  
    \[
        \langle c(\zeta_q, d) \rangle_{A_K} \subsetneq \mathcal{S}_K (\mathcal{G}_{n S})  
    \] 
    for every $d$ such that $q | d | n$. 
\end{theorem} 

\begin{proof}
    Choose $q,d$ such that $q \neq 1$ and $q | d | n$, and let $K$ be a splitting field for $\Phi_q(T)$. Choose some $\zeta_q \in K$ and consider $c = c (\zeta_q, d)$ as defined above. 
    First, observe that $c$ is a singular function: the sum of values of $c$ on the heads of the snake is $\Psi_d(\zeta_q)$, and as $q|d$ we have $\Psi_d (\zeta_q) = 0$. So $c \in \mathcal{S}_K$ by \cref{singular_headssumtozero}. 

    Second, we show that $\langle c \rangle_{A_K} \subsetneq \mathcal{S}_K$. By \cref{propersubid_oneonenotin}, it suffices to show that $[1,-1,0,...,0] \notin \langle c \rangle_{A_K}$. As before, $[1,-1,0,...,0] \in \langle c(\zeta_q, d) \rangle_{A_K}$ if and only if there exists $f \in A_K$ such that $f*c = [1,-1,0,...,0]$. We show that no such $f$ can exist. 
    Choose $f \in A_K$ and let $h:= f*c$. Then the value on head $i$ of $h$ is:  
    \begin{align} \label{hi}
        h_i =& \sum\limits_{j \in \Zmodn} f_{n-j+i} c_j \nonumber \\ 
        =& \sum\limits_{k \in [0,d-1]} f_{n-km+i} \zeta_q^k \ . \nonumber
    \end{align} 
     Assume for a contradiction that we have $h = [1,-1,0,...,0]$. Compare the values on the $0$th and $m$th heads: we must have $h_0 = 1$ and $h_m = 0$. Therefore:  
    \[
        h_0 = \sum\limits_{k \in [0,d-1]} f_{n-km} \zeta_q^k = 1 
    \] 
    and: 
    \begin{align*}
        h_m &= \sum\limits_{k \in [0,d-1]} f_{n-km+m} \zeta_q^k = 0 \\ 
            & = \left( \sum\limits_{k \in [0,d-1]} f_{n-(k-1)m} \zeta_q^{k-1} \right) \zeta_q \ . 
    \end{align*}

    Observe two things: first, consider the indices $n-(k-1)m$. Recall that $dm = n = 0$ in $\Zmodn$, and so:
    \begin{align*}
        k \in [0,d-1] \Rightarrow (k-1)m \in \{-m, 0, m, 2m, ..., (d-2)m\} = \{0,m,...,-2m,-m\} \ . 
    \end{align*}
    Second, as $q | d$, we have $\zeta_q^d=1$, and so for $\zeta_q^{k-1}$:  
    \begin{align*}
        \zeta_q^{k-1} & \in [\zeta_q^{-1}, \zeta_q^{d-2}] \\ 
        &\Rightarrow \zeta_q^r \in \{\zeta_q^{-1}, 1, \zeta_q ..., \zeta_q^{d} \zeta_q^{-2} \} = \{1, \zeta_q, ..., \zeta_q^{-2}, \zeta_q^{-1} \} \\ 
        &\Rightarrow \zeta_q^{r} \in [\zeta_q^{0}, \zeta_q^{d-1}] \ . 
    \end{align*} 
    That is, we can relabel $(k-1) \rightarrow j \in [0,d-1]$. Now: 
    \begin{align*}
        h_m &= \left( \sum\limits_{j \in [0,d-1]} f_{n-jm} \zeta_q^{j} \right) \zeta_q \\ 
            &= h_0 \zeta_q = 0 \ . 
    \end{align*} 
    However, we have $h_0 = 1$ and $\zeta_q \neq 0$; a contradiction. Therefore no such $f \in A_K$ exists, and so $\langle c(\zeta_q,d) \rangle_{A_K}$ is a proper subset ideal of $\mathcal{S}_K$. 
\end{proof} 

\begin{cor}
    If $K$ is a splitting field for $\Phi_q(T)$ for any $q \mid n$, then $A_K(\mathcal{G}_{nS})$ is not S-simple. 
\end{cor} 

When $n$ is even, we can have $q = 2$, with $\zeta_2 = -1$, thus giving us again the subset ideal of \cref{oneminusone_ideal} associated with $j = \frac{n}{2}$. 
In the prime-headed case $A_K(\mathcal{G}_{p S})$, there is only one choice for $q|d|p$ with $q \neq 1$; namely, $q=d=p$. Thus \cref{nheads_cycloideals} simplifies to: 
\begin{cor} \label{primeheads_cycloideals}
    Let $p$ be an odd prime. 
    If $K$ is a splitting field for the $p$th cyclotomic polynomial $\Phi_p(T)$, then $A_K (\mathcal{G}_{p S})$ has proper subset ideals given by:  
    \[
    \langle [1, \zeta_p, \zeta_p^2, ..., \zeta_p^{p-1}] \rangle_{A_K} \subsetneq \mathcal{S}_K (\mathcal{G}_{p S})
    \] 
    for each $\zeta_p$ primitive $p$th root of unity in $K$. 
\end{cor}

In fact, we can neatly enumerate all the elements of the ideal $\langle [1, \zeta_p, \zeta_p^2, ..., \zeta_p^{p-1}] \rangle_{A_K}$. Fix $\zeta_p \in K$. As in the proof of \cref{nheads_cycloideals}, choose $f \in A_K$ and let $h := f * [1,\zeta_p,...,\zeta_p^{p-1}]$. Then the value on head $i$ of $h$ is:  
    \begin{align*}
        h_i =& \sum\limits_{j \in \Zmodp} f_j \cdot \zeta_p^{p-j+i} \\ 
        =& \sum\limits_{j \in \Zmodn} f_j \cdot \zeta_p^{p-j} \cdot \zeta_p^i \\ 
        =& h_0 \zeta_p^i \ . 
    \end{align*}
    Therefore, $f * [1,\zeta_p,...,\zeta_p^{p-1}] = h_0 \cdot [1,\zeta_p,...,\zeta_p^{p-1}]$ where $h_0$ is a scalar that depends on the choice of $f$. And so: 
    \begin{equation} \label{cycloformideals_allscalarmultiples}
        \langle [1, \zeta_p, \zeta_p^2, ..., \zeta_p^{p-1}] \rangle_{A_K} = \{ k \cdot [1,\zeta_p,...,\zeta_p^{p-1}] : k \in K \} \ .
    \end{equation} 
Observe that each distinct root $\zeta_p$ gives a different subset ideal. Choose distinct roots $\zeta_p$ and $\zeta_p'$. It suffices to show that $[1, \zeta_p, \zeta_p^2, ..., \zeta_p^{p-1}] \notin  \langle [1, \zeta_p', \zeta_p'^{\ 2}, ..., \zeta_p'^{\ p-1}] \rangle_{A_K}$. 
Because $\zeta_p$ and $\zeta_p'$ are both primitive $p$th roots of unity, there exists precisely one $m$ such that $1 < m < p$ where $\zeta_p'^{\ m} = \zeta_p$. Therefore:  
 \[
    [1, \zeta_p, \zeta_p^2, ..., \zeta_p^{p-1}] = [1, \zeta_p'^{\ m}, \zeta_p'^{\ m+1}, ..., \zeta_p'^{\ p-m-1}] \ . 
 \] 
For this function to be an element of $\langle [1, \zeta_p', \zeta_p'^{\ 2}, ..., \zeta_p'^{\ p-1}] \rangle_{A_K}$, we would require some scalar $k$ such that $k \cdot [1, \zeta_p', \zeta_p'^{\ 2}, ..., \zeta_p'^{\ p-1}] = [1, \zeta_p'^{\ m}, \zeta_p'^{\ m+1}, ..., \zeta_p'^{\ p-m-1}] $. Comparing the first terms; we must have $k = 1$; comparing the final terms, we must have $k = \zeta_p'^{\ p-m}$. As $m < p$, we have that $\zeta_p'^{\ p-m} \neq 1$; a contradiction. 

The generalised $n$-headed version of \cref{cycloformideals_allscalarmultiples}, derived in a similar way as the proof of \cref{nheads_cycloideals}, is included for interest without proof: 
\begin{equation*}
    \langle c(\zeta_q, d) \rangle_{A_K} = \left\{ \sum\limits_{j \in [0,m-1]}  \left( k_j \cdot 1_{X_j}* c(\zeta_q, d) \right) : k_j \in K \right\} . 
\end{equation*}

\begin{remark}
    When $\mathrm{char}(K) = p$, we have that $1$ is a root of $\Phi_p(T) \in K[T]$. In the prime-headed case $A_K(\mathcal{G}_{p S})$, the cyclotomic-form subset ideals of \cref{primeheads_cycloideals} coincide with the $\langle [1,1,...,1] \rangle_{A_K}$ subset ideals of \cref{oneoneone_subideal}. 
    However, in the non-prime-headed case, when $\mathrm{char}(K) = p | n$, we have distinct subset ideals: 
    \[
         \langle [1,0,...,0, \ 1,0,...,0, \ ..., \ 1,0,...,0] \rangle_{A_K} \neq \langle [1,1,...,1] \rangle_{A_K} 
    \]
from \cref{nheads_cycloideals} and \cref{oneoneone_subideal} respectively. To see that these subset ideals are in fact different, recall from \cref{oneall_scalarmultiples} that the latter ideal contains only scalar multiples of $[1,1,...,1]$, and so cannot include the function $[1,0,...,0, \ ..., \ 1,0,...,0]$. 
\end{remark}   

\subsection{Summary} \label{summarycharts}
The three types of ideal that we have introduced are summarised here, along with the conditions under which they exist and are proper subset ideals of $A_K$. Let $K$ be any field and consider the algebra of the $n$-headed snake $A_K(\mathcal{G}_{n S})$ for any $n \ge 2$. 
\begin{table}[H] 
\centering 
\begin{tabular}{ ||c||c|c|c|| } 
 \hline
 & Result & Conditions & Proper Subset Ideal(s) \\ 
 \hline\hline 
 1 & \cref{oneoneone_subideal} & When $\mathrm{char}(K) | n$ & $\langle [1,1,...,1] \rangle_{A_K}$ \\ 
 \hline 
 2 & \cref{oneminusone_ideal} & When $n$ is not prime & $\langle 1_{X_0} - 1_{X_j} \rangle_{A_K} \ \forall \  j | n$  \\ 
 \hline 
 3 & \cref{nheads_cycloideals} & When $\zeta_q \in K$ for any $q | n$ & $\langle c(\zeta_q, d) \rangle_{A_K} \ \forall \ \zeta_q, \ \forall \ d$ such that $q|d|n$  \\ 
 \hline 
\end{tabular} 
\end{table}
For the prime-headed snake $\mathcal{G}_{p S}$, the subset ideal of \cref{oneminusone_ideal} does not occur, and the others simplify: 
\begin{table}[H] 
\centering 
\begin{tabular}{ ||c||c|c|c|| } 
 \hline
 & Result & Conditions (Prime heads) & Proper Subset Ideal(s) \\ 
 \hline\hline 
 1 & \cref{oneoneone_subideal} & When $\mathrm{char}(K) = p$ & $\langle [1,1,...,1] \rangle_{A_K}$ \\ 
 \hline 
 2 & \cref{oneminusone_ideal} & -- & -- \\ 
 \hline 
 3 & \cref{nheads_cycloideals} & When $\zeta_p \in K$ & $\langle [1, \zeta_p, \zeta_p^2, ..., \zeta_p^{p-1}] \rangle_{A_K} \ \forall \ \zeta_p$ \\ 
 \hline
\end{tabular} 
\end{table}

\subsection{Further Subset Ideals} \label{FurtherSubsetIdeals}

As the number of heads of the snake increases, so do the options for possible subset ideals. Consider the coefficients of the $n$ heads of $f \in A_K(\mathcal{G}_{n S})$ as $n$ different `variables' $f_j$. For an ideal generated by $s \in \mathcal{S}_K$, recall from \cref{licitvariations_EQUATION} that we have, without loss of generality, for a generating function:  
\[
    s = [1,s_1, ..., s_{n-2}, -(1+s_1+...+s_{n-2})] \ ,
\]
that is, up to $n-2$ nonzero `independent' variables. As in the proofs of \cref{oneoneone_subideal}, \cref{oneminusone_ideal}, and \cref{nheads_cycloideals}, we can calculate $f*s$ for an arbitrary $f \in A_K$. Setting the result equal to $[1,-1,0,..,0]$ gives a system of $n$ equations that can be solved for $f$ in terms of the fixed values $s_j$. We do so computationally; our code for the three, four, and five-headed cases is available on GitHub at \url{https://github.com/smarsland/SnakeCalculations} \cite{snakecode}.  

\begin{example} \label{3h_example} 
    In the three-headed snake, we consider most generally $s = [1,s_1, -(1+s_1)]$. In order for $[1,-1,0] \in \langle s \rangle_{A_K}$, we require there to exist an $f \in A_K$ such that: 
    \[
        f = \left[ f_0, f_0 - \frac{1}{s_1^2 + s_1 + 1}, f_0 + \frac{s_1}{s_1^2 + s_1 + 1} \right]
    \] 
    for any $f_0 \in K$. Such an $f$ exists if and only if the quantity $s_1^2 + s_1 + 1 \neq 0_K$. Therefore, proper subset ideals exist if and only if there is some $s_1 \in K$ which is a root of $\Phi_3(T)$, as in \cref{primeheads_cycloideals}. 
    Because we are able to compute this required $f$ for a general $\langle s \rangle_{A_K}$, we are able, in \cref{example:3h_ideals}, to completely characterise when $A_K (\mathcal{G}_{3 S})$ is and is not S-simple. 
\end{example} 
In this example we recover (a stronger statement of) \cref{primeheads_cycloideals}, which we elaborate on in \cref{Examples: Two and Three Headed Snakes} for the three-headed snake. With more heads on the snake, there are many more possibilities. 
Our computational method gives a `determinant' quantity $\Delta_s$ associated with each $s$ -- in the example above, $\Delta_s = s_1^2 + s_1 +1$. The solution for $f$ such that $f*s = [1,-1,0,...,0]$ is a scalar multiple of $\frac{1}{\Delta_s}$. Thus when $\Delta_s = 0_K$, such an $f$ cannot exist, and so $\langle s \rangle_{A_K}$ is a proper subset ideal. 
Because $s$ has up to $n-2$ `independent' variables, as the number of heads $n$ increases, so too does the difficulty of calculating $\Delta_s = 0_K$. In \cref{3h_example}, and in the three forms of subset ideal already detailed -- where the coefficients of the heads of $s$ are all dependent -- $\Delta_s$ is a univariate polynomial. In general, each additional head on the snake adds an additional potential variable to $\Delta_s$. 

For a general $n$-headed snake, and a general singular function $s \in \mathcal{S}_K(\mathcal{G}_{p S})$, we do not at this stage have an efficient method for calculating $\Delta_s$. But we can show that, even for the five-headed snake, there are many proper subset ideals not covered in \cref{ThreeConditionsSection}. 

Consider the singular functions of the five-headed snake that have support on precisely four heads: $\{ s = [1,s_1,s_2, -(s_1+s_2+1),0] : s_1,s_2 \neq 0 \text{ and } s_1+s_2 \neq -1 \}$. Using our code \cite{snakecode}, we solve $f*s = [1,-1,0,0,0]$ for $f$. Identifying the scalar factor of $f$ as $\frac{1}{\Delta_s}$, we have: 
\begin{align*}
    \Delta_s = \  &s_1^4 + s_1^3 + s_1^2 + s_1 + 1 \\ 
                &+ s_2 (2s_1^3 + 4s_1^2 + 6s_1 +3) \\ 
                &+ s_2^2 (4s_1^2 + 7s_1 + 4) \\ 
                &+ s_2^3 (3s_1 + 2) \\ 
                &+ s_2^4 \ .
\end{align*}
Setting $\Delta_s = 0_K$, fixing $s_1 \in K$, and solving for $s_2$, gives four potential values of $s_2$ per each $s_1 \in K$. 
Thus we obtain many options for $\langle [1,s_1,s_2, -(s_1+s_2+1),0]  \rangle_{A_K} \subsetneq \mathcal{S}_K$, depending on the field $K$. 
Note that these subset ideals are not covered by the three conditions already proved. As $n = 5$ is prime, there are no subgroup-associated subset ideals as in \cref{oneminusone_ideal}. With $s_4 = 0$, by \cref{oneall_scalarmultiples} and \cref{cycloformideals_allscalarmultiples} respectively we know that $s \notin \langle [1,1,1,1,1] \rangle_{A_K}$ and $s \notin \langle [1,\zeta_5, \zeta_5^2, \zeta_5^3, \zeta_5^4] \rangle_{A_K}$ for any $\zeta_5 \in K$. Interestingly, subset ideals of the form $\langle [1,s_1,s_2, -(s_1+s_2+1),0] \rangle_{A_K}$ can occur over fields that are not splitting fields of $\Phi_5(T)$; that is, fields over which we do not have the cyclotomic-form subset ideals of \cref{nheads_cycloideals}.

\begin{example}
    Let $K = \mathbb{F}_{59}$. We have that:  
    \[
        \langle [1,33,26,58,0] \rangle_{A_K}  \subsetneq \mathcal{S}_K (\mathcal{G}_{5 S}) 
    \]
    is a proper subset ideal in the five-headed snake that is distinct from the three types previously detailed. Note that $\mathbb{F}_{59}$ does not contain primitive fifth roots of unity $\zeta_5$ -- this can be confirmed computationally. 

    Using our code \cite{snakecode}, over $\mathbb{F}_{59}$ we find that:  
    \[
        f * [1,33,26,58,0] = [1,-1,0,0,0] \ \ \Leftrightarrow \ \ f = \frac{1}{49 \cdot 59}[k, 52+k, 52+k, k, 47+k ]
    \] 
    for any $k \in K$.   
    Thus $\Delta_s = 49 \cdot 59 = 0_K$, and so no such $f$ exists. Therefore $[1,-1,0,0,0] \notin \langle [1,33,26,58,0] \rangle_{A_K}$, and so $A_{\mathbb{F}_{59}} (\mathcal{G}_{5 S})$ is not S-simple. 
\end{example}

\begin{comment}
\begin{example}
    Let $K = \mathbb{F}_{31}$. Then, with $s = [1,5,6,-12,0]$, we have that:  
    \[
        \langle [1,5,6,-12,0] \rangle_{A_K}  \subsetneq \mathcal{S}_K (\mathcal{G}_{5 S}) 
    \]
    is a proper subset ideal in the five-headed snake that is distinct from the three types previously detailed.
    As $5$ is prime and $\mathrm{char}(\mathbb{F}_{31}) = 31 \nmid 5$, we have only to compare with the cyclotomic-form subset ideals of \cref{primeheads_cycloideals}. 

    Computationally, we find that:  
    \[
        f * [1,5,6,-12,0] = [1,-1,0,0,0] \ \ \Leftrightarrow \ \ f = \frac{1}{13051}[k, 227+k, -499+k, 1133+k, 53+k ]
    \] 
    for any $k \in K$.   
    In $\mathbb{F}_{31}$, this is:  
    \[
        f = \frac{1}{421 \cdot 31}[k,10+k, 28+k, 17+k, 22+k]
    \] 
    Thus $\Delta_s = 421 \cdot 31 = 0_K$, and so no such $f$ exists. Therefore $[1,-1,0,0,0] \notin \langle [1,5,6,-12,0] \rangle_{A_K}$. 

    This ideal is larger than the cyclotomic form ideals of \cref{primeheads_cycloideals}. The primitive fifth roots of unity in $\mathbb{F}_{31}$ are $\zeta_5 = 2,4,8$ and $16$. Calculation shows that three of the cyclotomic-form ideals (all of them save for $\zeta_5 = 8$) are contained in our new ideal, for example: 
    \[
        [1,2,4,8,16] = \frac{1}{18} [0,   21,    12,     20,     27 ] * [1,5,6,-12,0] 
    \] 
    Further separate calculations for each of the four roots $\zeta_5$ can confirm that $s = [1,5,6,-12,0]$ is not contained in any of the four cyclotomic-form ideals. 
\end{example} 
\end{comment}

Even for five-head-supported functions $s$, the equation $\Delta_s$ becomes impractically long; and for an $n$-headed snake, $\Delta_s$ contains terms up to a power of $n-1$. 
Further research in this direction would be computational in focus. For this paper, it suffices to note that our three conditions under which $A_K(\mathcal{G}_{n S})$ is not S-simple -- $\mathrm{char}(K) | n$; $n$ being non-prime; and $K$ being a splitting field for $\Phi_n(T)$ -- are far from exhaustive.

\section{Example: the Three-Headed Snake} \label{Examples: Two and Three Headed Snakes}

We conclude with a full overview of the most straightforward example where subset ideals exist, the three-headed snake, in which we \emph{can} completely characterise all possible proper subset ideals of $A_K(\mathcal{G}_{3 S})$. 

The two-headed case, dealt with in \cref{2h_alwaysssimple}, is trivial:  $A_K(\mathcal{G}_{2S})$ is always S-simple, regardless of the field $K$. 
The three-headed case is more interesting. As $3$ is a prime, there are never subset ideals of the form of \cref{oneminusone_ideal}. Unlike in the five-headed snake, we show that the two other results of \cref{ThreeConditionsSection} are the only other possible proper subset ideals of $A_K(\mathcal{G}_{3 S})$. 

\begin{prop} \label{example:3h_ideals} 
    For the three-headed snake, $A_K (\mathcal{S}_{3 S})$ is non-S-simple if and only if $K$ is a splitting field of $\Phi_3(T)$.  
\end{prop}

\begin{proof} 
The proof follows directly from \cref{3h_example} above, which tells us computationally that any proper subideals of $\mathcal{S}_K (\mathcal{G}_{3 S})$ are of the form $\langle [1, b , -(b+1)] \rangle_{A_K}$ for $b \in K$ such that $\Phi_3(b) = 0_K$. 
Therefore $A_K(\mathcal{G}_{3 S})$ is not S-simple if and only if the root(s) $\zeta_3$ are elements of $K$; that is, if and only if $\Phi_3(T)$ splits over $K$. 
\end{proof} 

Note that if $\mathrm{char}(K) = 3$, then $K$ is automatically a splitting field for $\Phi_3(T)$; both roots $\zeta_3$ are equal to $1$ -- see \cref{3h_repeatedrootonlyone} -- and so the functions $[1,\zeta_3, \zeta_3^2], [1,\zeta_3', \zeta_3'^2]$ and $[1,1,1]$ coincide. In this case there is only one proper subset ideal, $\langle [1,1,1] \rangle_{A_K}$ of \cref{oneoneone_subideal}. If $K$ is otherwise a splitting field for $\Phi_3(T)$, the two values of $\zeta_3$ are distinct, and there are two proper subset ideals, of the cyclotomic form of \cref{primeheads_cycloideals}. 

\subsection{Splitting fields of $\Phi_3(T)$} 

The splitting fields of cyclotomic polynomials are known as \emph{cyclotomic fields}, and are widely used. For an overview, see for example \cite{DummitandFoote}. 
As the three-headed snake is of particular interest as a `complete' example,  we conclude by providing a classification of the splitting fields of $\Phi_{3}(T)$, using elementary number-theoretic arguments. 

\subsubsection{Case 1 -- Splitting fields of characteristic zero}
First, we consider the case where $\mathrm{char}(K) = 0$: this implies that $K$ contains a subfield isomorphic\footnote{We assume that such an isomorphism is fixed once for all, and we will write $\Q \subseteq K$ accordingly.} to $\Q$. Fix an algebraic closure $\overline{K}$ of $K$ that contains $\overline{\Q}$, and let $i, \sqrt{3} \in \overline{\Q} \subseteq \overline{K}$ be fixed square roots of $-1$ and $3$ respectively. We then have the following result.

\begin{prop}
The field $K$ is a splitting field of $\Phi_{3}(T)$ if and only if it contains $\frac{i\sqrt{3}}{2}$ -- that is, if $K$ is a field extension of $\Q\left[\frac{i\sqrt{3}}{2}\right]$.

In this case, the two singular ideals of $A_K(\GGG_{3S})$ are: 
    \[
        \left\langle \left[1,\frac{-1 + i\sqrt{3}}{2} , \frac{-1 - i\sqrt{3}}{2} \right] \right\rangle_{A_K} \text{and } \ \left\langle \left[1, \frac{-1 - i\sqrt{3}}{2}, \frac{-1 + i\sqrt{3}}{2}\right] \right\rangle_{A_K}  
    \]
\end{prop}

\begin{proof}
In $\overline{K}[T]$ we have:
\[ \begin{array}{rcl} \displaystyle \left(T - \frac{-1 + i\sqrt{3}}{2} \right)\left(T - \frac{-1-i\sqrt{3}}{2}\right) & = &  \displaystyle T^{2} + \frac{1 - i\sqrt{3} + 1 + i\sqrt{3}}{2} T + \frac{(1 - i\sqrt{3})(1 + i\sqrt{3})}{4} \\
& = & T^{2} + T + 1 \\
& = & \Phi_{3}(T) \ , 
\end{array}\]
and so the roots of $\Phi_{3}(T)$ in $\overline{K}$ are $\displaystyle \frac{-1 + i\sqrt{3}}{2} \text{ and } \frac{-1 - i\sqrt{3}}{2}$.
\end{proof}

\subsubsection{Case 2 -- Splitting fields of positive characteristic} 
Assume now that $K$ is a splitting field of $\Phi_{3}(T)$ of positive characteristic $p$, and let $\F_{p} = \Zmodp$ be the prime subfield of $K$. The following results are standard. 

\begin{lemma} \label{3h_repeatedrootonlyone}
    Let $K$ be a splitting field for $\Phi_3(T)$. Then there is a repeated root of $\Phi_3(T)$ if and only if $\mathrm{char}(K) = 3$, in which case this root is $1_K$; otherwise the two roots $\zeta_3$ and $\zeta_3'$ are distinct and not equal to $1$. 
\end{lemma}

\begin{proof}
    First, assume $\mathrm{char}(K) = 3$. Then $\F_3$ is the prime subfield of $K$, and 
    \[\displaystyle (T - 1)^{2} = T^{2} - 2T + 1 = T^{2} + T + 1 \ \text{ in $\F_3[T] \subseteq K[T]$ . } \]
    Thus $1$ is the only root (with multiplicity $2$) of $\Phi_{3}(T)$ in $K$. Observe that if $\zeta_3 = 1$, then $\Phi_3(1) = 3 = 0_K$ also implies that $\mathrm{char}(K) = 3$. 

    Second, assume that $\Phi_3(T)$ has a repeated root. 
    Let $b$ be a root of $\Phi_3(T)$; that is, $b^2 + b + 1 = 0_K$. The other root is $-(b+1)$: 
    \[
        (-(b+1))^2 - (b+1) + 1 = b^2 + 2b + 1 -b = 0 \ . 
    \] 
    When $b = -(b+1)$ we have $2b = -1$; that is, $b = -(2^{-1})$ in $K$. For $\Phi_3(-(2^{-1})) = \frac{3}{4} = 0_K$ we must have $\mathrm{char}(K) = 3$. 
\end{proof}

\begin{lemma} \noindent  
\begin{enumerate}
\label{3h_Deltabfactors0or1mod3}
    \item \label{it1:galois} 
    Let $b \in \mathbb{Z}$. 
    All factors of $b^{2}+b+1$ are congruent to either $0 \ \mathrm{mod}\  3$ or $1 \ \mathrm{mod} \ 3$. 
    \item \label{it2:galois} The field $\F_{p}$ is a splitting field for $T^{2}+T+1$ if and only if $p = 3$ or $p \equiv 1 \ \mathrm{mod} \ 3$. 
    \end{enumerate}
\end{lemma}

\begin{proof}
    First we prove \eqref{it1:galois}. If $q$ is a prime factor of $b^{2}+b+1$, then $q$ is odd as $b^{2}$ and $b$ have the same parity, hence $b^{2} + b + 1$ is odd. Moreover, $q$ must also divide
    \[ 4 (b^{2}+b+1) = (2b+1)^{2} + 3, \] 
    so $-3$ must be a quadratic residue mod $q$. Hence the corresponding Legendre symbol $\left(\frac{-3}{q}\right)$ must be equal to $1$. The multiplicativity of the Legendre symbol and  the Quadratic Reciprocity Law applied to the pair of odd primes $(3,q)$ then imply that $q$ must satisfy: 
    \begin{equation} \label{LRQ}
\biggl(\frac{-1}{q}\biggr)\biggl(\frac{q}{3}\biggr) = \Bigl(-1\Bigr)_.^{\frac{q-1}{2}} 
    \end{equation}
Depending on whether $-1$ is a square mod $q$, i.e. whether we have $q \equiv 1 \ \mathrm{mod} \ 4$, we have the following dichotomy: 
\begin{itemize}
\item If $-1$ is a square mod $q$, i.e. if $q \equiv 1 \ \mathrm{mod} \ 4$, then \cref{LRQ} becomes:
\[\displaystyle\left(\frac{q}{3}\right) = 1 \ , \text{ i.e. $q$ is a square mod $3$.} \]
As $q$ is a prime number, it cannot be divisible by $3$ unless it equals $3$ (in which case it is a square mod $3$); as the only squares mod $3$ are $0$ and $1$, we obtain that either $q = 3$, or $q$ is congruent to $1$ mod $3$.
\item If $-1$ is not a square mod $q$, i.e. if $q \equiv -1 \ \mathrm{mod} \ 4$, then \cref{LRQ} becomes:
\[\displaystyle -\left(\frac{q}{3}\right) = -1, \ \text{ i.e.} \left(\frac{q}{3}\right) = 1 ,\]
and we can use again the previous argument to prove that either $q = 3$, or $q$ is congruent to $1$ mod $3$.
\end{itemize} 

Now we prove \eqref{it2:galois}. Suppose that $\mathbb{F}_p$ is a splitting field for $T^{2}+T+1$. Then there exists an integer $b \in \Z$ such that $b^{2}+b+1 \equiv 0 \ \mathrm{mod} \ p$. This means that $p$ is a prime factor of $b^{2}+b+1$ in $\mathbb{Z}$, and \eqref{it1:galois} above ensures that either $p = 3$ or $p \equiv 1 \ \mathrm{mod} \ 3$, as claimed. 

First, assume that  $p = 3$. Then $\mathrm{char}(\mathbb{F}_3) = 3$, so by \cref{3h_repeatedrootonlyone} we have that $\mathbb{F}_3$ is a splitting field of $\Phi_{3}(T)$.
Second, assume that $p \equiv 1 \ \mathrm{mod} \ 3$. Then $p$ is a quadratic residue mod $3$, so the multiplicativity of the Legendre symbol and the Quadratic Reciprocity law ensure that: 
    \begin{equation} \label{CSLRQ}
    \displaystyle \biggl(\frac{-3}{p}\biggr) = \biggl(\frac{-1}{p}\biggr)\biggl(\frac{3}{p}\biggr) = \biggl(\frac{-1}{p}\biggr)\biggl(\frac{p}{3}\biggr) \Bigl(-1 \Bigr)^{\frac{p-1}{2}} = \Bigl(-1 \Bigr)^{\frac{p-1}{2}} \biggl(\frac{-1}{p}\biggr) . 
    \end{equation}
    Since $-1$ is a quadratic residue mod $p$ if and only if $p \equiv 1 \ \mathrm{mod} \ 4$, which is itself equivalent to $(-1)^{\frac{p-1}{2}} = 1$, we obtain that the right-hand side of \cref{CSLRQ} is equal to $1$, which proves that $-3$ is a quadratic residue mod $p$.

    Now: choose an integer $n$ such that $p$ divides $n^{2} + 3$. Replacing $n$ by $n + p$ if necessary, we can assume without loss of generality that $n$ is odd, and write $n = 2b + 1$ with $b \in \Z$. Thus, we have:
    \[\displaystyle (2b+1)^{2} \equiv -3 \ \mathrm{mod} \ p \ , \text{ i.e. } \ 4b^{2} + 4b +4 \equiv 0 \ \mathrm{mod} \ p \ .\]
    As $p$ is odd (hence coprime to $4$), this proves that such an integer $b$ satisfies: 
    \[ \displaystyle b^{2} + b + 1 \equiv 0 \ \mathrm{mod} \ p \ . \]
    The image $b \in \F_{p}$ is thus a root of $\Phi_{3}(T) = T^{2} + T + 1$ in $\F_{p}$. Thus $\F_{p}$ is a splitting field of $\Phi_{3}(T)$. 
\end{proof}

\begin{cor} \label{3h:CSnonsimple}
    If the characteristic of $K$ is either $3$ or a prime integer congruent to $1$ modulo $3$, then $A_{K}(\mathcal{G}_{3 S})$ is not S-simple. 
\end{cor}

\begin{proof}
Let $p >0$ be the characteristic of $K$, and $\F_{p} \subseteq K$ be the prime subfield of $K$. Our assumption on $p$ ensures by \cref{3h_Deltabfactors0or1mod3} that $\F_{p}$, and therefore $K$, are splitting fields of $\Phi_{3}(T)$. The conclusion follows from \cref{primeheads_cycloideals}. 
\end{proof}

\begin{remark}
\cref{3h:CSnonsimple} provides a sufficient condition but not a necessary one. Indeed, for example, if $K$ is an algebraic closure of $\F_{2} = \Z/2\Z$, then $K$ is a field of positive characteristic $2$: it does not satisfy the assumption of \cref{3h:CSnonsimple}. Nevertheless, as $K$ is algebraically closed, \cref{primeheads_cycloideals} shows that $A_{K}(\GGG_{3S})$ is not S-simple. 
\end{remark}

\bibliographystyle{plain} 
\bibliography{Bibliography} 

\end{document}